\documentclass[11pt]{article}
\usepackage{amssymb}
\newtheorem{theorem}{Theorem}[section]

\newcommand\Hecke{{\cal{H}}}
\newcommand\integers{{\mathbb{Z}}}
\newcommand\rationals{{\mathbb{Q}}}
\newcommand\complexes{{\mathbb{C}}}

\begin{document}

%\headline={\hss \tenrm --\  \folio\ -- \hss}
%\footline={\hfill}

%{version:  \number\month{/}\number\day{/}\number\year}

\title 
{ Subrings of the asymptotic Hecke algebra of type $H_4$} 

\author{
Dean Alvis\\
Department of Mathematical Sciences\\
Indiana~University~South~Bend\\
South~Bend, IN, 46634, USA\\
}

% \begin{abstract}
% Abstract. \ 
% The structure of subring $J^{\Gamma \cap \Gamma^{-1}}$ 
% of the asymptotic Hecke algebra is 
% described for $\Gamma$ a left cell of the Coxeter
% group of type $H_4$.  A small set of generators 
% over $\integers$ is produced.  
% The subalgebras spanned by a subset of the basis 
% $\left\{t_x\right\}_{x\in \Gamma\cap\Gamma^{-1}}$
% are determined.
% \end{abstract}

\maketitle

%%%%%%%%%%%%%%%%%%%%%%%%%%%%%%%%%%%%%%%%%%%%%%%%%%%%%%
%%%%%%%%%%%%%%%%%%%%%%%%%%%%%%%%%%%%%%%%%%%%%%%%%%%%%%
%%   Section 1

\section{Introduction}

Let $W$ be a finite Coxeter
group with set of distinguished generators $S$,
length function $\ell:w \mapsto \ell(w)$,
and Bruhat order $\le$.
Let $J$ be the {\it asymptotic Hecke algebra}
of $W$, as defined by Lusztig in 
\cite[Section~2]{LusztigAffineTwo} (see also 
\cite[Chapter~18]{LusztigTwoParam}).  As an
additive group, 
$J$ is a free abelian group with basis 
$\left(t_w\right)_{w\in W}$ indexed by $W$.  
The multiplication 
operation of $J$ is given by 
\[
t_x t_y
=
\sum_{z \in W}
  \gamma_{x,y,z^{-1}} t_z,
\]
where the structure constants
$\gamma_{x,y,z^{-1}} \in \integers$
are described in the next section.
It is known that $J$ is an associative ring
with identity.  
Moreover,
if $\Gamma$ is a left cell of $W$, then
\[
J^{\Gamma \cap \Gamma^{-1}}
=
\sum_{x \in \Gamma \cap \Gamma^{-1}} 
 \integers \, t_x
\]
is a $\integers$-subalgebra of $J$.
We denote this ring $J(\Gamma)$.

Fokko~du~Cloux
has computed $\gamma_{x,y,z^{-1}}$
for all $x,y,z \in W = W(H_4)$.
In fact, 
du~Cloux has determined all of the coefficients,
not just the leading coefficients,
of the structure constants $h_{x,y,z}$ of
the Hecke algrebra: see \cite{Cloux}.
By du~Cloux's calculations, the coefficients
of the $h_{x,y,z}$ are nonnegative integers. 
Since the same is known for the
Kazhdan-Lusztig polynomials, 
results of Lusztig (\cite[Chapter 15]{LusztigTwoParam})
show that all of the conjectures of 
\cite[Chapter 14]{LusztigTwoParam},
 hold in type $H_4$.
In particular, each left cell $\Gamma$
of $W(H_4)$ contains a unique element
of $\cal{D}$, the set of distinguished
involutions.  Moreover, if 
$e \in \Gamma \cap \cal{D}$,
then $t_e$ is the identity element of
$J(\Gamma)$.

In the current investigation, the
structure constants $\gamma_{x,y,z^{-1}}$
are needed only
for $x,y,z \in \Gamma \cap \Gamma^{-1}$,
$\Gamma$ a left cell of $W=W(H_4)$. 
These constants
were calculated by computer using
the algorithm described in Section~2, 
which differs from that used by du~Cloux.
There are a total of 206 left cells in
type $H_4$ (\cite{AlvisHFour}).  
For $\Gamma$ a
left cell of $W$, the associated  
$W$-graph gives rise to a corresponding
$\rationals W$-module as in \cite{KazhdanLusztig}, 
which will be denoted 
$M(\Gamma)$.  
We say a bijection 
$\pi : \Gamma_1 \cap \Gamma_1^{-1} \rightarrow
       \Gamma_2 \cap \Gamma_2^{-1}$
is a 
{\it permutation isomorphism}
from $J(\Gamma_1)$ to $J(\Gamma_2)$ if
\[
\gamma_{x,y,z^{-1}} = \gamma_{\pi(x),\pi(y),\pi(z)^{-1}}
\]
for all
$x,y,z \in \Gamma_1 \cap \Gamma_1^{-1}$.
A computer search of the matrices
of structure constants reveals the following
result.  (The author knows of no a priori 
proof of this result.)  

\begin{theorem}
Suppose $(W,S)$ is of type $H_4$ and
$\Gamma_1$, $\Gamma_2$ are left cells of $W$
such that the corresponding modules 
$M(\Gamma_1)$, $M(\Gamma_2)$ are isomorphic.
Then there is a unique permutation
isomorphism 
from
$J(\Gamma_1)$ to
$J(\Gamma_2)$.
\label{permisothm}
\end{theorem}

For a left cell $\Gamma$ not in the largest  
two-sided cell $A$, we have 
$\vert \Gamma \cap \Gamma^{-1} \vert = 1$ or $2$,
and hence $J(\Gamma)$ is easily described (see Section~3).  
Thus the interesting cases are those
for which $\Gamma \subseteq A$.  There are
three isomorphism classes of modules
$M(\Gamma)$, $\Gamma \subseteq A$, represented
by $\Gamma = A_1$, $A_9$, $A_{19}$ (in the notation
of \cite{AlvisHFour}).
Sections~4--6 describe the rings $J(\Gamma)$
for these cases.  A CAS program was used to find
a set of generators over $\integers$, the 
characteristic polynomials for the left multiplication
operators $(t_x)_L : J(\Gamma) \rightarrow J(\Gamma)$, 
and the subsets of 
$\left\{t_x\right\}_{x \in\Gamma\cap\Gamma^{-1}}$
spanning subalgebras of $J(\Gamma)$.

The author is indebted to Victor Ostrik for suggesting
this problem, and to George Lusztig and Victor Ostrik
for several helpful communications.

%%%%%%%%%%%%%%%%%%%%%%%%%%%%%%%%%%%%%%%%%%%%%%%%%%%%%%
%%%%%%%%%%%%%%%%%%%%%%%%%%%%%%%%%%%%%%%%%%%%%%%%%%%%%%
%%   Section 2

\section{The computation of the structure constants}

Let $W$, $S$ be as in the previous section.  
Let $\Hecke$ be the corresponding Hecke algebra
over $A = \integers[q^{1/2},q^{-1/2}]$,
$q$ an indeterminate, with standard basis
$\left(T_w\right)_{w\in W}$ satisfying
\begin{equation}
T_s T_w
=
\left\{
\begin{array}{cl}
T_{sw} & \mbox{if $sw > w$} \\
T_{sw}^{\strut} + (q^{1/2} - q^{-1/2})T_w 
         & \mbox{if $sw < w$} 
\end{array}
\right.
\label{stdbasiseqn}
\end{equation}
for $s\in S$, $w \in W$.
(This notation of \cite{LusztigTwoParam} differs
slightly from that in \cite{KazhdanLusztig}.)
The semilinear involution $a \mapsto \overline{a}$
of $\Hecke$ is given by $\overline{q^{1/2}}=q^{-1/2}$,
$\overline{T_w} = T_{w^{-1}}^{-1}$.
The basis
$\left( c_w \right)_{w\in W}$ for
$\Hecke$ (denoted 
$\left(C_w^{\prime}\right)_{w\in W}$ 
in \cite{KazhdanLusztig}) satisfies
\[
c_w = \sum_{y \in W} p_{y,w} T_y
\]
where $p_{y,w} \in q^{-1/2}\integers[q^{-1/2}]$
when $y<w$, $p_{w,w}=1$, $p_{y,w}=0$ when $y \not\le w$,
and
$\overline{c_w} = c_w$.  

For $x,y,z\in W$,
define $f^\prime_{x,y,z} \in A$ by
\[
T_x T_y
=
\sum_{z \in W}
f^\prime_{x,y,z} c_z.
\]
Then $\gamma_{x,y,z^{-1}}$ is determined
by
\begin{equation}
f^\prime_{x,y,z}
=
\gamma_{x,y,z^{-1}} q^{a(z)/2}
+
\mbox{lower degree terms}
\label{fprimeeqn}
\end{equation}
(\cite[13.6(d)]{LusztigTwoParam}), where
$a(z)$ is a nonnegative integer depending
only on the two-sided cell containing $z$
(see below).
Now, if $f_{x,y,z} \in A$ are given by
\[
T_x T_y = \sum_{z \in W} f_{x,y,z} T_z,
\]
then
\[
f^\prime_{x,y,z}
=
\sum_{w \in W} p^{\prime}_{z,w} f_{x,y,w}
\]
by \cite[13.1(b)]{LusztigTwoParam},
where
$\left[p^{\prime}_{z,w}\right]$ is the inverse matrix
of
$\left[p_{z,w}\right]$. 
Further, if $W$ is finite then
\[
p_{z,w}^{\prime} 
= \varepsilon_z \varepsilon_w p_{w_0w, w_0z}
\]
by \cite[11.4]{LusztigTwoParam}, 
where $w_0$ is the longest element of $W$ and 
$\varepsilon_x = (-1)^{\ell(x)}$.
Put
$q_x^{1/2} = (q^{1/2})^{\ell(x)}$,  
and let
\[
P_{x,y} =  q_x^{-1/2} q_y^{1/2} p_{x,y},
\]
so $P_{x,y}$
is the Kazhdan-Lusztig polynomial for $x$, $y$.
Define 
\[
F_{x,y,z} = \  q_x^{1/2} q_y^{1/2} q_z^{-1/2} f_{x,y,z}.
\]
Then 
\begin{eqnarray*}
f_{x,y,z}^\prime
&=& 
\sum_{w \in W}
\left( \varepsilon_z \varepsilon_{w}
   q_{w}^{-1/2} q_z^{1/2} P_{w_0w,w_0z}
\right)
\left(
 q_x^{-1/2} q_y^{-1/2} q_{w}^{1/2} F_{x,y,w}
\right) \\
&=& 
q_x^{-1/2} q_y^{-1/2} q_z^{1/2}
\sum_{w \in W}
\varepsilon_z \varepsilon_{w}
     P_{w_0w,w_0z} F_{x,y,w}. 
\end{eqnarray*}
Therefore
formula (\ref{fprimeeqn}) is equivalent to
\begin{equation}
\sum_{w \in W}
\varepsilon_z \varepsilon_{w}
     P_{w_0w,w_0z}\, F_{x,y,w}
=
\gamma_{x,y,z^{-1}} q^{\left(a(z)-\ell(x)-\ell(y)+\ell(z)\right)/2}
+
\mbox{lower degree terms}.
\label{neweqn}
\end{equation}

To find the structure constants 
$\gamma_{x,y,z^{-1}}$
for $x,y,z \in\Gamma \cap \Gamma^{-1}$ in type $H_4$, 
the polynomials
$F_{x,y,w}$ were evaluated by computer for a
fixed $x \in \Gamma\cap\Gamma^{-1}$ 
and all $y \in \Gamma\cap\Gamma^{-1}$, $w \in W$, 
using a
straightforward calculation based on 
(\ref{stdbasiseqn}).
The leading term of the sum on the left
side of (\ref{neweqn}) was then found for
$y,z \in \Gamma\cap \Gamma^{-1}$, using the
Kazhdan-Lusztig polynomials computed in the
course of determining the left cells 
in \cite{AlvisHFour}.  
Varying $x$ 
over $\Gamma\cap \Gamma^{-1}$ produced the
value of the $a$-function on $\Gamma$: 
if $\delta(x,y,z)$ denotes the degree of the
left side of (\ref{neweqn}) and
\[
\mu = \max
\left\{
2 \delta(x,y,z) - \left(\ell(x) + \ell(y) - \ell(z)\right)
\mid
x,y,z \in \Gamma \cap\Gamma^{-1}
\right\},
\]
then $a(x) = \mu$ for $x \in \Gamma\cap\Gamma^{-1}$.
Once the value $a(x)$ had been found,
the structure constants $\gamma_{x,y,z^{-1}}$
were then determined using (\ref{neweqn}).  
The results of these calculations are summarized
in the next sections. 

%%%%%%%%%%%%%%%%%%%%%%%%%%%%%%%%%%%%%%%%%%%%%%%%%%%%%%
%%%%%%%%%%%%%%%%%%%%%%%%%%%%%%%%%%%%%%%%%%%%%%%%%%%%%%
%%   Section 3

\section{Small left cells}

For the remainder of this paper 
$(W,S)$ is of type $H_4$.
The notations of \cite{AlvisHFour} are used 
for the left and two-sided cells.  
(There is a typographic error in \cite{AlvisHFour}: 
the left cell $A_{12}$ is equal to $A_{11}d$,
not $A_{10}d$.) 
In particular, the two-sided cells of $W$
are $A$, $B$, $C$, $D$, $E$, $F$, $G$, with
$A$ the ``big'' cell.

If $\Gamma$ is contained in one
of the two-sided cells $B$, $C$ or $G$, then
$\vert \Gamma \cap \Gamma^{-1} \vert = 1$.
In this case $J(\Gamma) = \integers t_e$
where $e \in \Gamma\cap \Gamma^{-1}$,
and $t_e^2 = t_e$.

Now suppose $\Gamma$ is contained in one
of the two-sided cells $D$, $E$, or $F$.
In this case
$\vert \Gamma \cap \Gamma^{-1} \vert = 2$, 
and
$\Gamma \cap \Gamma^{-1}
= \left\{e,s \right\}$
where $e$ is the distinguished involution
 and $s$ is the other involution in $\Gamma$.
Then
\[
J(\Gamma)
=
\integers t_e \oplus \integers t_s,
\]
with identity element $t_e$.
Moreover, the calculations described in 
Section~2 show 
\[
t_s^2
=
\left\{
\begin{array}{cl}
t_e & \mbox{if $\Gamma \subseteq D$,} \\
t_e + t_s^{\strut} & \mbox{if $\Gamma \subseteq E$ or
       $\Gamma \subseteq F$.} \\
\end{array}
\right.
\]

From these results and the structure of the
modules $M(\Gamma)$ given in \cite{AlvisHFour},
 Theorem~\ref{permisothm} holds
for left cells $\Gamma$ not contained in $A$.

%%%%%%%%%%%%%%%%%%%%%%%%%%%%%%%%%%%%%%%%%%%%%%%%%%%%%%
%%%%%%%%%%%%%%%%%%%%%%%%%%%%%%%%%%%%%%%%%%%%%%%%%%%%%%
%%   Section 4

\section{The case $\Gamma = A_1$}

It remains only to consider the left cells 
$\Gamma$ such that
$\Gamma \subseteq A$.  
Suppose $\Gamma$ is the left cell $A_1$, so
$\vert \Gamma \vert = 326$ and
$\Gamma \cap \Gamma^{-1} = 14$ (\cite{AlvisHFour}).
Let
$S=\left\{a,b,c,d\right\}$, where
$(ab)^3=(bc)^3=(cd)^5=(ac)^2=(ad)^2=(bd)^2=1$.  
The elements $x_1$, \dots, $x_{14}$ of
$\Gamma \cap \Gamma^{-1}$ are indexed 
according to the
list of reduced expressions given in 
Table~\ref{a1elements}.

\begin{table}
\caption{The elements of $\Gamma \cap \Gamma^{-1}$, $\Gamma = A_1$.}
\[
\vbox{
\offinterlineskip
\halign{
\strut $#$ & \quad $#$ \cr
j &  x_j \hfil \cr
\noalign{\hrule} \cr
1 & abcaba \cr
2 & abcdabcaba \cr
3 & abcdabcdabcaba \cr
4 & abcdabcdabcdabcaba \cr
5 & bcdabcdabcdbcdabcaba \cr
6 & abcdabcdabcdabcdabcaba \cr
7 & abcdbcdabcdabcdbcdabcaba \cr
8 & bcdabcdabcdbcdabcdabcaba \cr
9 & abcdabcdabcdabcdabcdabcaba \cr
10 & abcdbcdabcdabcdbcdabcdabcaba \cr
11 & abcdabcdabcdabcdabcdabcdabcaba \cr
12 & abcdabcdabcdabcdabcdabcdabcdabcaba \cr
13 & abcdabcdabcdabcdabcdabcdabcdabcdabcaba \cr
14 & abcdabcdabcdabcdabcdabcdabcdabcdabcdabcaba \cr
}}
\]
\label{a1elements}
\end{table}

The structure constants
$\gamma_{x,y,z^{-1}}$ are described below
by giving matrices
$M_1$, \dots, $M_{14}$, where
for fixed $j$, $M_j$ is the
matrix $\left[\gamma_{x_j,y,z^{-1}}\right]$,
with $y$ and $z$ varying over
$\Gamma \cap \Gamma^{-1}$ in the order given in 
Table~\ref{a1elements}.
Note that $M_j$ is the transpose of the
left multiplication operator
$(x_j)_L : J(\Gamma) \rightarrow J(\Gamma)$.
To save space, only a set of generators is
given explicitly, and the other matrices
are then described in terms of those
generators.

The calculations described in
Section~2 yield $M_1 = I$, the identity matrix,
so $x_1$ is the distinguished involution of $A_1$.
Also, 
\[
M_{2} = 
{%\matrixfont
\left(
\begin{array}{cccccccccccccc}
0&1&0&0&0&0&0&0&0&0&0&0&0&0\cr
1&0&1&0&0&0&0&0&0&0&0&0&0&0\cr
0&1&1&1&0&0&0&0&0&0&0&0&0&0\cr
0&0&1&1&1&1&0&0&0&0&0&0&0&0\cr
0&0&0&1&0&0&1&0&0&0&0&0&0&0\cr
0&0&0&1&0&1&0&1&1&0&0&0&0&0\cr
0&0&0&0&1&0&0&0&1&0&0&0&0&0\cr
0&0&0&0&0&1&0&0&0&1&0&0&0&0\cr
0&0&0&0&0&1&1&0&1&0&1&0&0&0\cr
0&0&0&0&0&0&0&1&0&0&1&0&0&0\cr
0&0&0&0&0&0&0&0&1&1&1&1&0&0\cr
0&0&0&0&0&0&0&0&0&0&1&1&1&0\cr
0&0&0&0&0&0&0&0&0&0&0&1&0&1\cr
0&0&0&0&0&0&0&0&0&0&0&0&1&0\cr
\end{array}
\right)
}, 
\]
\[
M_{4} = 
{%\matrixfont
\left(
\begin{array}{cccccccccccccc}
0&0&0&1&0&0&0&0&0&0&0&0&0&0\cr
0&0&1&1&1&1&0&0&0&0&0&0&0&0\cr
0&1&2&3&1&2&1&1&1&0&0&0&0&0\cr
1&1&3&4&2&4&1&1&3&1&1&0&0&0\cr
0&1&1&2&0&2&1&1&1&0&1&0&0&0\cr
0&1&2&4&2&4&2&2&4&1&3&1&0&0\cr
0&0&1&1&1&2&0&1&2&1&1&1&0&0\cr
0&0&1&1&1&2&1&0&2&1&1&1&0&0\cr
0&0&1&3&1&4&2&2&4&2&4&2&1&0\cr
0&0&0&1&0&1&1&1&2&0&2&1&1&0\cr
0&0&0&1&1&3&1&1&4&2&4&3&1&1\cr
0&0&0&0&0&1&1&1&2&1&3&2&1&0\cr
0&0&0&0&0&0&0&0&1&1&1&1&0&0\cr
0&0&0&0&0&0&0&0&0&0&1&0&0&0\cr
\end{array}
\right)
},
\]
and
\[
M_{6} = 
{%\matrixfont
\left(
\begin{array}{cccccccccccccc}
0&0&0&0&0&1&0&0&0&0&0&0&0&0\cr
0&0&0&1&0&1&1&0&1&0&0&0&0&0\cr
0&0&1&2&1&3&1&1&2&1&1&0&0&0\cr
0&1&2&4&2&4&2&2&4&1&3&1&0&0\cr
0&0&1&2&1&2&0&1&2&1&1&1&0&0\cr
1&1&3&4&2&5&2&2&5&2&4&2&1&0\cr
0&1&1&2&0&2&1&1&2&1&2&1&0&0\cr
0&0&1&2&1&2&1&1&2&0&2&1&1&0\cr
0&1&2&4&2&5&2&2&5&2&4&3&1&1\cr
0&0&1&1&1&2&1&0&2&1&2&1&0&0\cr
0&0&1&3&1&4&2&2&4&2&4&2&1&0\cr
0&0&0&1&1&2&1&1&3&1&2&1&0&0\cr
0&0&0&0&0&1&0&1&1&0&1&0&0&0\cr
0&0&0&0&0&0&0&0&1&0&0&0&0&0\cr
\end{array}
\right)
}.
\]
Moreover,
\[
M_{3} =  -I +M_{2}^{2},
\quad
M_{5} =  I -M_{4} -M_{6} -M_{2}^{2} +M_{2}M_{4},
\]
\[
M_{7} =  I -M_{2} -2M_{4} -M_{2}^{2} -M_{2}M_{4} -M_{6}M_{2} +M_{2}^{2}M_{4},
\]
\[
M_{8} =  I -M_{2} -2M_{4} -M_{2}^{2} -M_{2}M_{4} -M_{2}M_{6} +M_{2}^{2}M_{4},
\]
\[
M_{9} =  -I +M_{2} +M_{4} -M_{6} +M_{2}^{2} +M_{2}M_{4} +M_{2}M_{6} +M_{6}M_{2} -M_{2}^{2}M_{4},
\]
\begin{eqnarray*}
M_{10} 
&= &
 I -2M_{2} -3M_{4} +2M_{6} -2M_{2}^{2} -4M_{2}M_{4} -2M_{2}M_{6} \\
& & \qquad {}+2M_{4}^{2} +M_{4}M_{6} -2M_{6}M_{2} 
			+3M_{2}^{2}M_{4} -M_{2}M_{4}^{2},
\end{eqnarray*}
\[
M_{11} =  2M_{4} -M_{6} +M_{2}M_{4} -M_{4}^{2} -M_{4}M_{6} -2M_{2}^{2}M_{4} +M_{2}M_{4}^{2},
\]
\[
M_{12} =  -I +M_{2} -M_{4} +3M_{6} +2M_{2}^{2} -M_{2}M_{4} +M_{4}^{2} +3M_{4}M_{6} +3M_{2}^{2}M_{4} -2M_{2}M_{4}^{2},
\]
\begin{eqnarray*}
M_{13} 
&= & I +M_{4} -4M_{6} -2M_{2}^{2} +3M_{2}M_{4} +M_{2}M_{6} \\
& &\qquad {}-2M_{4}^{2} -4M_{4}M_{6} +M_{6}M_{2} +M_{6}^{2} -3M_{2}^{2}M_{4} +2M_{2}M_{4}^{2}, 
\end{eqnarray*}
and
\begin{eqnarray*}
M_{14} 
& = &
 2I -M_{2} +M_{4} -6M_{6} -2M_{2}^{2} +4M_{2}M_{4} +M_{2}M_{6} \\
& & \qquad {}-2M_{4}^{2} -6M_{4}M_{6} +M_{6}M_{2} -2M_{6}^{2} -5M_{2}^{2}M_{4} +3M_{2}M_{4}^{2} +M_{2}M_{6}^{2}.
\end{eqnarray*}
This completes the description of the multiplication structure
constants when $\Gamma = A_1$.  
After similar calculations were carried out
for the left cells $A_2$,\dots, $A_8$,
a computer search established the
corresponding cases of Theorem~\ref{permisothm}.  
In addition, we have the following.

\begin{theorem}
Suppose $\Gamma = A_1$ and
$X$ is a nonempty subset of $\Gamma\cap \Gamma^{-1}$.
Then 
$\sum_{x \in X}\integers t_x$ is a subring (with 1) of
$J(\Gamma)$ if and only if
$X$ is one of the sets 
$\left\{x_1\right\}$,
$\left\{x_1,x_{14}\right\}$,
$\Gamma\cap \Gamma^{-1}$.
\end{theorem}

The characteristic polynomials of
the matrices $M_j$ appear in Table~\ref{a1polys}.

\begin{table}
\caption{Characteristic polynomials for $\Gamma=A_1$.}
\[
\vbox{
\offinterlineskip
\halign{
\strut $#$ & \quad $#$ \hfil \cr
j & \det(u I - M_j) \cr
\noalign{\hrule}1 & (-1+u)^{14}\cr
2 & (-2+u)(-1+u)u(1+u)(-2+u^{2})^{2}(1-3u+u^{2})(-4-2u+u^{2})(-1+u+u^{2})\cr
3 & (-3+u)(-1+u)^{4}u^{2}(1+u)(5-10u+u^{2})(-5-5u+u^{2})(-1-u+u^{2})\cr
4 & (-1+u)^{3}(1+u)^{5}(1-18u+u^{2})(1-7u+u^{2})(1+3u+u^{2})\cr
5 & u^{2}(1+u)^{2}(-2+u^{2})^{2}(-4-8u+u^{2})(5-5u+u^{2})(1-3u+u^{2})\cr
6 & (-1+u)u(1+u)(2+u)(-2+u^{2})^{2}(-4-22u+u^{2})(1-3u+u^{2})(-1-u+u^{2})\cr
7 & u^{4}(1+u)^{2}(3+u^{2})^{2}(5-10u+u^{2})(-4+2u+u^{2})\cr
8 & u^{4}(1+u)^{2}(3+u^{2})^{2}(5-10u+u^{2})(-4+2u+u^{2})\cr
9 & (-1+u)^{2}u(2+u)(-2+u^{2})^{2}(-4-22u+u^{2})(-1-u+u^{2})(1+3u+u^{2})\cr
10 & (-1+u)u^{2}(1+u)(-2+u^{2})^{2}(-4-8u+u^{2})(5-5u+u^{2})(1+3u+u^{2})\cr
11 & (-1+u)^{6}(1+u)^{2}(1-18u+u^{2})(1+3u+u^{2})(1+7u+u^{2})\cr
12 & (-3+u)u^{2}(1+u)^{5}(5-10u+u^{2})(-1-u+u^{2})(-5+5u+u^{2})\cr
13 & (-2+u)(-1+u)^{2}u(-2+u^{2})^{2}(-4-2u+u^{2})(-1+u+u^{2})(1+3u+u^{2})\cr
14 & (-1+u)^{7}(1+u)^{7}\cr
}}
\]
\label{a1polys}
\end{table}

%%%%%%%%%%%%%%%%%%%%%%%%%%%%%%%%%%%%%%%%%%%%%%%%%%%%%%
%%%%%%%%%%%%%%%%%%%%%%%%%%%%%%%%%%%%%%%%%%%%%%%%%%%%%%
%%   Section 5

\section{The case $\Gamma = A_9$}

In this section $\Gamma$ is the left cell $A_9$
of $W$, so
$\vert \Gamma \vert = 392$ and
$\Gamma \cap \Gamma^{-1} = 18$.
We use a adapt the notation of the previous section
to this case.   
The elements $x_1$, \dots, $x_{18}$ of
$\Gamma \cap \Gamma^{-1}$ are indexed as in 
Table~\ref{a9elements}.

\begin{table}
\caption{The elements of $\Gamma \cap \Gamma^{-1}$, $\Gamma = A_9$.}
\[
\vbox{
\offinterlineskip
\halign{
\strut $#$ & \quad $#$ \cr
j & x_j \cr
\noalign{\hrule}
1 & cdcdca \cr
2 & dcdabcdbcdca \cr
3 & cdcdabcdabcdca \cr
4 & dcdabcdbcdbcdbca \cr
5 & cdcdabcdabcdabcdca \cr
6 & cdcdbcdabcdabcdbcdca \cr
7 & dcdabcdbcdabcdabcdca \cr
8 & cdcdabcdabcdabcdabcdca \cr
9 & dcdabcdbcdabcdabcdabcdca \cr
10 & cdabcdbcdbcdabcdabcdbcdbca \cr
11 & cdcdabcdabcdabcdabcdabcdca \cr
12 & cdcdabcdabcdbcdabcdabcdbcdca \cr
13 & dcdabcdbcdabcdabcdabcdabcdca \cr
14 & cdcdabcdabcdabcdabcdabcdabcdca \cr
15 & cdcdabcdabcdabcdabcdabcdabcdabcdca \cr
16 & dcdabcdbcdabcdabcdbcdabcdabcdbcdca \cr
17 & dcdabcdbcdabcdabcdabcdabcdabcdabcdca \cr
18 & cdcdabcdabcdabcdabcdabcdabcdabcdabcdabcdca \cr
}}
\]
\label{a9elements}
\end{table}

Let $M_j$ be the $18\times18$-matrix
$\left[\gamma_{x_j,y,z^{-1}}\right]$.
Then 
$M_1 = I$,
\[
M_{2} = 
{%\matrixfont
\left(
\begin{array}{cccccccccccccccccc}
0&1&0&0&0&0&0&0&0&0&0&0&0&0&0&0&0&0\cr
1&1&1&1&1&0&0&0&0&0&0&0&0&0&0&0&0&0\cr
0&1&1&0&1&0&1&0&0&0&0&0&0&0&0&0&0&0\cr
0&1&0&0&0&1&0&0&0&0&0&0&0&0&0&0&0&0\cr
0&1&1&0&2&1&1&1&1&0&0&0&0&0&0&0&0&0\cr
0&0&0&1&1&1&0&0&1&1&0&0&0&0&0&0&0&0\cr
0&0&1&0&1&0&0&1&0&0&1&0&0&0&0&0&0&0\cr
0&0&0&0&1&0&1&1&1&0&1&0&1&0&0&0&0&0\cr
0&0&0&0&1&1&0&1&2&0&1&1&0&1&0&0&0&0\cr
0&0&0&0&0&1&0&0&0&0&0&1&0&0&0&0&0&0\cr
0&0&0&0&0&0&1&1&1&0&1&0&1&1&0&0&0&0\cr
0&0&0&0&0&0&0&0&1&1&0&1&0&1&0&1&0&0\cr
0&0&0&0&0&0&0&1&0&0&1&0&0&1&1&0&0&0\cr
0&0&0&0&0&0&0&0&1&0&1&1&1&2&1&0&1&0\cr
0&0&0&0&0&0&0&0&0&0&0&0&1&1&1&0&1&0\cr
0&0&0&0&0&0&0&0&0&0&0&1&0&0&0&0&1&0\cr
0&0&0&0&0&0&0&0&0&0&0&0&0&1&1&1&1&1\cr
0&0&0&0&0&0&0&0&0&0&0&0&0&0&0&0&1&0\cr
\end{array}
\right)
},
\]
\[
M_{3} = 
{%\matrixfont
\left(
\begin{array}{cccccccccccccccccc}
0&0&1&0&0&0&0&0&0&0&0&0&0&0&0&0&0&0\cr
0&1&1&0&1&1&0&0&0&0&0&0&0&0&0&0&0&0\cr
1&1&1&0&1&0&0&1&0&0&0&0&0&0&0&0&0&0\cr
0&0&0&0&1&0&0&0&0&0&0&0&0&0&0&0&0&0\cr
0&1&1&1&2&1&1&1&1&0&1&0&0&0&0&0&0&0\cr
0&1&0&0&1&1&0&0&1&0&0&1&0&0&0&0&0&0\cr
0&0&0&0&1&0&1&1&1&0&0&0&0&0&0&0&0&0\cr
0&0&1&0&1&0&1&1&1&0&1&0&0&1&0&0&0&0\cr
0&0&0&0&1&1&1&1&2&1&1&1&1&1&0&0&0&0\cr
0&0&0&0&0&0&0&0&1&0&0&0&0&0&0&0&0&0\cr
0&0&0&0&1&0&0&1&1&0&1&0&1&1&1&0&0&0\cr
0&0&0&0&0&1&0&0&1&0&0&1&0&1&0&0&1&0\cr
0&0&0&0&0&0&0&0&1&0&1&0&1&1&0&0&0&0\cr
0&0&0&0&0&0&0&1&1&0&1&1&1&2&1&1&1&0\cr
0&0&0&0&0&0&0&0&0&0&1&0&0&1&1&0&1&1\cr
0&0&0&0&0&0&0&0&0&0&0&0&0&1&0&0&0&0\cr
0&0&0&0&0&0&0&0&0&0&0&1&0&1&1&0&1&0\cr
0&0&0&0&0&0&0&0&0&0&0&0&0&0&1&0&0&0\cr
\end{array}
\right)
},
\]
and
\[
M_{4} = 
{%\matrixfont
\left(
\begin{array}{cccccccccccccccccc}
0&0&0&1&0&0&0&0&0&0&0&0&0&0&0&0&0&0\cr
0&1&0&0&0&0&1&0&0&0&0&0&0&0&0&0&0&0\cr
0&0&0&0&1&0&0&0&0&0&0&0&0&0&0&0&0&0\cr
1&0&0&1&0&0&0&0&0&1&0&0&0&0&0&0&0&0\cr
0&0&1&0&1&0&0&0&1&0&0&0&0&0&0&0&0&0\cr
0&0&0&0&0&1&0&1&0&0&0&0&0&0&0&0&0&0\cr
0&1&0&0&0&0&1&0&0&0&0&0&1&0&0&0&0&0\cr
0&0&0&0&0&1&0&1&0&0&1&0&0&0&0&0&0&0\cr
0&0&0&0&1&0&0&0&1&0&0&0&0&1&0&0&0&0\cr
0&0&0&1&0&0&0&0&0&1&0&0&0&0&0&1&0&0\cr
0&0&0&0&0&0&0&1&0&0&1&1&0&0&0&0&0&0\cr
0&0&0&0&0&0&0&0&0&0&1&1&0&0&0&0&0&0\cr
0&0&0&0&0&0&1&0&0&0&0&0&1&0&0&0&1&0\cr
0&0&0&0&0&0&0&0&1&0&0&0&0&1&1&0&0&0\cr
0&0&0&0&0&0&0&0&0&0&0&0&0&1&0&0&0&0\cr
0&0&0&0&0&0&0&0&0&1&0&0&0&0&0&1&0&1\cr
0&0&0&0&0&0&0&0&0&0&0&0&1&0&0&0&1&0\cr
0&0&0&0&0&0&0&0&0&0&0&0&0&0&0&1&0&0\cr
\end{array}
\right)
}
\]
by the calculations described in 
Section~2.
Further, 
\[
M_{5} =  -I -M_{2} -M_{3} -M_{4} +M_{2}^{2},
\quad
M_{6} =  I +M_{4} -M_{2}^{2} +M_{2}M_{3},
\]
\[
M_{7} =  -M_{2} +M_{2}M_{4}, \quad
M_{8} =  M_{4} -M_{2}^{2} +M_{3}^{2},
\]
\[
M_{9} =  I +M_{2} +M_{3} -M_{2}^{2} -2M_{2}M_{3} -2M_{2}M_{4} -M_{3}^{2} +M_{2}^{3},
\]
\[
M_{10} =  -I -M_{4} +M_{4}^{2}, \quad
M_{11} =  I +2M_{2} +M_{2}M_{4} -M_{3}^{2} -M_{2}^{3} +M_{2}^{2}M_{3},
\]
\[
M_{12} =  -2I -2M_{2} -M_{4} +2M_{2}^{2} +M_{2}M_{3} +2M_{2}M_{4} +M_{3}^{2} -M_{2}^{3} -M_{2}^{2}M_{3} +M_{2}M_{3}^{2},
\]
\[
M_{13} =  -2M_{2}M_{4} +M_{2}M_{4}^{2},
\]
\[
M_{14} =  -2I -2M_{2} -2M_{3} -M_{4} +2M_{2}^{2} +2M_{2}M_{3} +M_{3}^{2} -2M_{2}^{2}M_{3} +M_{3}^{3},
\]
\begin{eqnarray*}
M_{15} 
& =  &
3I +3M_{2} +5M_{3} +3M_{4} -M_{2}^{2} +4M_{2}M_{3} +6M_{2}M_{4} +M_{3}^{2} -M_{4}^{2} \\
& &\qquad {}-M_{2}^{3} +2M_{2}^{2}M_{3} -3M_{2}M_{3}^{2} -2M_{2}M_{4}^{2} -4M_{3}^{3} +M_{2}^{2}M_{3}^{2},
\end{eqnarray*}
\[
M_{16} =  I -M_{4} -2M_{4}^{2} +M_{4}^{3},
\]
\begin{eqnarray*}
M_{17} 
& =  &
5I +6M_{2} +2M_{3} +4M_{4} -4M_{2}^{2} -3M_{2}M_{3} +2M_{2}M_{4} -M_{4}^{2} +6M_{2}^{2}M_{3} \\
& &\qquad {}-3M_{2}M_{3}^{2} -M_{2}M_{4}^{2} -2M_{3}^{3} -M_{2}^{2}M_{3}^{2} +M_{2}M_{3}^{3},
\end{eqnarray*}
and
\[
M_{18} =  3M_{4} -3M_{4}^{3} +M_{4}^{4}.
\]
This completes the description of the structure constants
for the left cell $\Gamma = A_9$.
As before, after the structure constants for
$A_{10}$,\dots,$A_{18}$ were also computed,
the relevant cases of Theorem~\ref{permisothm} were verified by a
computer search.
We also have the following.

\begin{theorem}
Suppose $\Gamma = A_9$  and
$X$ is a nonempty subset of $\Gamma\cap \Gamma^{-1}$.
Then 
$\sum_{x \in X}\integers t_x$ is a subring (with 1) of
$J(\Gamma)$ if and only if
$X$ is one of the sets
\[
\left\{x_1\right\}, \quad
\left\{x_1,x_{18}\right\}, \quad
\left\{x_1,x_4,x_{10},x_{16},x_{18}\right\}, \quad
\Gamma\cap \Gamma^{-1}.
\]
\end{theorem}

The characteristic polynomials of
the matrices $M_j$ appear in
Table~\ref{a9polys}.

\begin{table}
\caption{Characteristic polynomials for $\Gamma=A_9$.}
\[
\vbox{
\offinterlineskip
\halign{
\strut $#$ & \quad $#$ \hfil \cr
j & \det(u I - M_j) \cr
\noalign{\hrule}1 & (-1+u)^{18}\cr
2 & (-1+u)^{2}u^{2}(-5-5u+u^{2})(-1-4u+u^{2})(-5-u+u^{2})^{2}(-1-u+u^{2})^{2}(-1+u+u^{2})\cr
3 & (-1+u)(1+u)^{3}(1-7u+u^{2})(-1-4u+u^{2})(1-3u+u^{2})^{2}(-1-u+u^{2})^{2}(-1+u+u^{2})\cr
4 & (1+u)^{2}(1-3u+u^{2})^{4}(-1-u+u^{2})^{4}\cr
5 & (-1+u)^{7}(1+u)^{5}(1-18u+u^{2})(1-7u+u^{2})(1+3u+u^{2})\cr
6 & u^{2}(1+u)^{2}(1+u^{2})^{2}(5-10u+u^{2})(1-3u+u^{2})(-1-u+u^{2})(5+2u+u^{2})^{2}\cr
7 & u^{2}(1+u)^{2}(1+u^{2})^{2}(5-10u+u^{2})(1-3u+u^{2})(-1-u+u^{2})(5+2u+u^{2})^{2}\cr
8 & u^{2}(-5-15u+u^{2})(1-3u+u^{2})(-1-u+u^{2})^{4}(-5+u+u^{2})^{2}\cr
9 & (-1+u)u^{8}(1+u)(-4-22u+u^{2})(-4+2u+u^{2})^{3}\cr
10 & (-1+u)^{2}u^{8}(-4-2u+u^{2})^{4}\cr
11 & u^{2}(-5-15u+u^{2})(1-3u+u^{2})(-5+u+u^{2})^{2}(-1+u+u^{2})^{4}\cr
12 & (-1+u)^{2}u^{2}(1+u^{2})^{2}(5-10u+u^{2})(-1-u+u^{2})(5+2u+u^{2})^{2}(1+3u+u^{2})\cr
13 & (-1+u)^{2}u^{2}(1+u^{2})^{2}(5-10u+u^{2})(-1-u+u^{2})(5+2u+u^{2})^{2}(1+3u+u^{2})\cr
14 & (-1+u)^{9}(1+u)^{3}(1-18u+u^{2})(1+3u+u^{2})(1+7u+u^{2})\cr
15 & (-1+u)(1+u)^{3}(1-7u+u^{2})(1-3u+u^{2})^{2}(-1-u+u^{2})(-1+u+u^{2})^{2}(-1+4u+u^{2})\cr
16 & (1+u)^{2}(1-3u+u^{2})^{4}(-1+u+u^{2})^{4}\cr
17 & (-1+u)^{2}u^{2}(-5-5u+u^{2})(-5-u+u^{2})^{2}(-1-u+u^{2})(-1+u+u^{2})^{2}(-1+4u+u^{2})\cr
18 & (-1+u)^{10}(1+u)^{8}\cr
}}
\]
\label{a9polys}
\end{table}

%%%%%%%%%%%%%%%%%%%%%%%%%%%%%%%%%%%%%%%%%%%%%%%%%%%%%%
%%%%%%%%%%%%%%%%%%%%%%%%%%%%%%%%%%%%%%%%%%%%%%%%%%%%%%
%%   Section 6

\section{The case $\Gamma = A_{19}$}

Suppose $\Gamma$ is the left cell $A_{19}$.
Thus 
$\vert \Gamma \vert = 436$ and
$\Gamma \cap \Gamma^{-1} = 24$.
A notation similar to that in the previous
two sections is used for the elements of 
$\Gamma \cap \Gamma^{-1}$ and the matrices
of structure constants.
Table~\ref{a19elements} lists the 
elements $x_1$, \dots, $x_{24}$ of
$\Gamma \cap \Gamma^{-1}$.

\begin{table}
\caption{The elements of $\Gamma \cap \Gamma^{-1}$, $\Gamma = A_{19}$.}
\[
\vbox{
\offinterlineskip
\halign{
\strut $#$ & \quad $#$ \cr
j & x_j \cr
\noalign{\hrule}
1 & dcdbcdbcdc \cr
2 & dcdbcdbcdbcdbc \cr
3 & cdcdbcdabcdabcdc \cr
4 & cdbcdbcdbcdabcdabcdc \cr
5 & cdcdbcdabcdabcdabcdc \cr
6 & cdcdbcdabcdabcdbcdbc \cr
7 & dcdbcdbcdabcdabcdabcdc \cr
8 & cdbcdbcdbcdabcdabcdbcdbc \cr
9 & cdcdbcdabcdabcdabcdabcdc \cr
10 & cdcdbcdabcdbcdabcdabcdbcdc \cr
11 & dcdbcdabcdabcdbcdabcdabcdc \cr
12 & dcdbcdbcdabcdabcdabcdabcdc \cr
13 & cdcdbcdabcdabcdabcdabcdabcdc \cr
14 & dcdbcdbcdabcdabcdabcdabcdabcdc \cr
15 & dcdbcdabcdabcdbcdabcdabcdbcdbc \cr
16 & cdcdbcdabcdabcdabcdabcdabcdabcdc \cr
17 & dcdbcdbcdabcdabcdbcdabcdabcdbcdc \cr
18 & dcdbcdbcdabcdabcdabcdabcdabcdabcdc \cr
19 & cdcdbcdabcdabcdabcdabcdabcdabcdabcdc \cr
20 & cdcdbcdabcdbcdabcdabcdbcdabcdabcdabcdc \cr
21 & dcdbcdbcdabcdabcdabcdabcdabcdabcdabcdc \cr
22 & cdcdbcdabcdabcdabcdabcdabcdabcdabcdabcdc \cr
23 & cdcdbcdabcdbcdabcdabcdbcdabcdabcdabcdabcdc \cr
24 & cdcdbcdabcdabcdabcdabcdabcdabcdabcdabcdabcdc \cr
}}
\]
\label{a19elements}
\end{table}

Let $M_j$ be the $24\times 24$ matrix
$\left[\gamma_{x_j,y,z^{-1}}\right]$.
Then 
$M_2 = I$,
\[
M_{4} = 
{%\matrixfont
\left(
\begin{array}{cccccccccccccccccccccccc}
0&0&0&1&0&0&0&1&0&0&0&0&0&0&0&0&0&0&0&0&0&0&0&0\cr
0&0&0&1&0&0&0&0&0&0&0&0&0&0&0&0&0&0&0&0&0&0&0&0\cr
1&0&0&1&0&0&1&0&0&0&0&0&0&1&0&0&0&0&0&0&0&0&0&0\cr
0&0&0&0&0&0&1&0&0&0&0&0&0&0&0&0&0&0&0&0&0&0&0&0\cr
0&0&1&0&0&0&1&0&0&0&0&0&1&0&0&0&0&0&0&0&0&0&0&0\cr
0&1&0&0&1&0&0&1&0&0&0&0&0&0&0&0&1&0&0&0&0&0&0&0\cr
0&0&1&0&1&1&1&0&0&0&0&1&1&0&0&0&0&0&0&0&0&0&0&0\cr
0&0&0&0&1&0&0&0&0&0&0&0&0&0&0&0&0&0&0&0&0&0&0&0\cr
0&0&0&0&1&0&0&0&0&0&0&1&0&0&0&0&0&1&0&0&0&0&0&0\cr
0&0&0&0&0&0&1&0&0&0&0&0&1&0&0&1&0&0&0&0&0&0&0&0\cr
0&0&0&1&0&0&0&0&1&1&0&0&0&1&0&0&0&0&0&0&1&0&0&0\cr
0&0&0&0&1&0&0&0&0&1&1&1&0&0&0&0&0&1&0&0&0&0&0&0\cr
0&0&0&0&0&0&1&0&1&1&1&0&1&0&1&1&0&0&0&0&0&0&0&0\cr
0&0&0&0&0&0&0&0&1&1&0&0&0&0&0&0&0&0&0&0&0&0&0&0\cr
0&0&0&0&0&0&0&1&0&1&0&0&0&0&0&0&1&0&0&0&0&1&0&0\cr
0&0&0&0&0&0&0&0&0&0&0&1&1&0&0&1&0&1&1&1&0&0&0&0\cr
0&0&0&0&0&0&0&0&0&1&0&0&0&0&0&0&0&0&0&0&0&0&0&0\cr
0&0&0&0&0&0&0&0&0&0&0&0&1&0&0&1&0&0&1&0&0&0&0&0\cr
0&0&0&0&0&0&0&0&0&0&0&0&0&1&0&1&0&0&0&0&1&0&0&1\cr
0&0&0&0&0&0&0&0&0&0&0&0&0&0&0&0&1&1&0&0&0&1&1&0\cr
0&0&0&0&0&0&0&0&0&0&0&0&0&0&0&1&0&0&0&0&0&0&0&0\cr
0&0&0&0&0&0&0&0&0&0&0&0&0&0&0&0&0&1&0&0&0&0&0&0\cr
0&0&0&0&0&0&0&0&0&0&0&0&0&0&0&0&0&0&0&0&1&0&0&0\cr
0&0&0&0&0&0&0&0&0&0&0&0&0&0&0&0&0&0&0&0&1&1&0&0\cr
\end{array}
\right)
},
\]
\[
M_{8} = 
{%\matrixfont
\left(
\begin{array}{cccccccccccccccccccccccc}
0&0&0&1&0&0&0&0&0&0&0&0&0&0&0&0&0&0&0&0&0&0&0&0\cr
0&0&0&0&0&0&0&1&0&0&0&0&0&0&0&0&0&0&0&0&0&0&0&0\cr
0&0&0&0&0&0&1&0&0&0&0&0&0&0&0&0&0&0&0&0&0&0&0&0\cr
1&0&0&1&0&0&0&0&0&0&0&0&0&1&0&0&0&0&0&0&0&0&0&0\cr
0&0&0&0&1&1&0&0&0&0&0&1&0&0&0&0&0&0&0&0&0&0&0&0\cr
0&0&0&0&1&0&0&0&0&0&0&0&0&0&0&0&0&0&0&0&0&0&0&0\cr
0&0&1&0&0&0&1&0&0&0&0&0&1&0&0&0&0&0&0&0&0&0&0&0\cr
0&1&0&0&0&0&0&1&0&0&0&0&0&0&0&0&1&0&0&0&0&0&0&0\cr
0&0&0&0&0&0&0&0&0&1&1&0&0&0&0&0&0&0&0&0&0&0&0&0\cr
0&0&0&0&0&0&0&0&1&1&1&0&0&0&1&0&0&0&0&0&0&0&0&0\cr
0&0&0&0&0&0&0&0&1&1&0&0&0&0&0&0&0&0&0&0&0&0&0&0\cr
0&0&0&0&1&0&0&0&0&0&0&1&0&0&0&0&0&1&0&0&0&0&0&0\cr
0&0&0&0&0&0&1&0&0&0&0&0&1&0&0&1&0&0&0&0&0&0&0&0\cr
0&0&0&1&0&0&0&0&0&0&0&0&0&1&0&0&0&0&0&0&1&0&0&0\cr
0&0&0&0&0&0&0&0&0&1&0&0&0&0&0&0&0&0&0&0&0&0&0&0\cr
0&0&0&0&0&0&0&0&0&0&0&0&1&0&0&1&0&0&1&0&0&0&0&0\cr
0&0&0&0&0&0&0&1&0&0&0&0&0&0&0&0&1&0&0&0&0&1&0&0\cr
0&0&0&0&0&0&0&0&0&0&0&1&0&0&0&0&0&1&0&1&0&0&0&0\cr
0&0&0&0&0&0&0&0&0&0&0&0&0&0&0&1&0&0&0&0&0&0&0&0\cr
0&0&0&0&0&0&0&0&0&0&0&0&0&0&0&0&0&1&0&0&0&0&0&0\cr
0&0&0&0&0&0&0&0&0&0&0&0&0&1&0&0&0&0&0&0&1&0&0&1\cr
0&0&0&0&0&0&0&0&0&0&0&0&0&0&0&0&1&0&0&0&0&1&1&0\cr
0&0&0&0&0&0&0&0&0&0&0&0&0&0&0&0&0&0&0&0&0&1&0&0\cr
0&0&0&0&0&0&0&0&0&0&0&0&0&0&0&0&0&0&0&0&1&0&0&0\cr
\end{array}
\right)
},
\]
and
\[
M_{21} = 
{%\matrixfont
\left(
\begin{array}{cccccccccccccccccccccccc}
0&0&0&0&0&0&0&0&0&0&0&0&0&0&0&0&0&0&0&0&1&1&0&0\cr
0&0&0&0&0&0&0&0&0&0&0&0&0&0&0&0&0&0&0&0&1&0&0&0\cr
0&0&0&0&0&0&0&0&0&0&0&0&0&1&0&1&0&0&0&0&1&0&0&1\cr
0&0&0&0&0&0&0&0&0&0&0&0&0&0&0&1&0&0&0&0&0&0&0&0\cr
0&0&0&0&0&0&0&0&0&0&0&0&1&0&0&1&0&0&1&0&0&0&0&0\cr
0&0&0&0&0&0&0&0&0&0&0&0&0&0&0&0&1&1&0&0&0&1&1&0\cr
0&0&0&0&0&0&0&0&0&0&0&1&1&0&0&1&0&1&1&1&0&0&0&0\cr
0&0&0&0&0&0&0&0&0&0&0&0&0&0&0&0&0&1&0&0&0&0&0&0\cr
0&0&0&0&1&0&0&0&0&0&0&1&0&0&0&0&0&1&0&0&0&0&0&0\cr
0&0&0&0&0&0&1&0&0&0&0&0&1&0&0&1&0&0&0&0&0&0&0&0\cr
0&0&0&1&0&0&0&0&1&1&0&0&0&1&0&0&0&0&0&0&1&0&0&0\cr
0&0&0&0&1&0&0&0&0&1&1&1&0&0&0&0&0&1&0&0&0&0&0&0\cr
0&0&0&0&0&0&1&0&1&1&1&0&1&0&1&1&0&0&0&0&0&0&0&0\cr
0&0&0&0&0&0&0&0&1&1&0&0&0&0&0&0&0&0&0&0&0&0&0&0\cr
0&0&0&0&0&0&0&1&0&1&0&0&0&0&0&0&1&0&0&0&0&1&0&0\cr
0&0&1&0&1&1&1&0&0&0&0&1&1&0&0&0&0&0&0&0&0&0&0&0\cr
0&0&0&0&0&0&0&0&0&1&0&0&0&0&0&0&0&0&0&0&0&0&0&0\cr
0&0&1&0&0&0&1&0&0&0&0&0&1&0&0&0&0&0&0&0&0&0&0&0\cr
1&0&0&1&0&0&1&0&0&0&0&0&0&1&0&0&0&0&0&0&0&0&0&0\cr
0&1&0&0&1&0&0&1&0&0&0&0&0&0&0&0&1&0&0&0&0&0&0&0\cr
0&0&0&0&0&0&1&0&0&0&0&0&0&0&0&0&0&0&0&0&0&0&0&0\cr
0&0&0&0&1&0&0&0&0&0&0&0&0&0&0&0&0&0&0&0&0&0&0&0\cr
0&0&0&1&0&0&0&0&0&0&0&0&0&0&0&0&0&0&0&0&0&0&0&0\cr
0&0&0&1&0&0&0&1&0&0&0&0&0&0&0&0&0&0&0&0&0&0&0&0\cr
\end{array}
\right)
}
\]
Also,
\[
M_{1} =  M_{4} +M_{21} +2M_{4}M_{8} -M_{4}M_{8}^{2},
\quad
M_{3} =  M_{4}^{2} +M_{4}M_{21} +2M_{4}^{2}M_{8} -M_{4}^{2}M_{8}^{2},
\]
\[
M_{5} =  M_{8}M_{4}, \quad
M_{6} =  M_{8}M_{4} +M_{8}M_{21} +2M_{4}^{3} -2M_{4}^{2}M_{8} -M_{4}^{3}M_{8} +M_{4}^{2}M_{8}^{2},
\]
\[
M_{7} =  M_{4}^{2}, \quad
M_{9} =  -M_{8} -M_{4}^{2} +M_{8}M_{4} +M_{4}M_{8}M_{4} -M_{8}^{2}M_{4},
\]
\[
M_{10} =  -M_{4} -M_{8}M_{4} +M_{8}^{2}M_{4},
\]
\[
M_{11} =  M_{4} +M_{8}M_{4} -M_{8}^{2} +M_{4}^{3} -2M_{4}^{2}M_{8} -M_{8}^{2}M_{4} +M_{4}^{4} -M_{4}^{3}M_{8} -M_{4}^{2}M_{8}M_{4} +M_{4}^{2}M_{8}^{2},
\]
\[
M_{12} =  -2M_{8}M_{4} -M_{8}M_{21} -M_{4}^{3} +M_{4}^{2}M_{8} +M_{4}^{3}M_{8} -M_{4}^{2}M_{8}^{2},
\]
\[
M_{13} =  -2M_{4}^{2} -M_{4}M_{21} -M_{4}^{2}M_{8} +M_{4}^{2}M_{8}^{2},
\quad
M_{14} =  -2M_{4} -M_{21} -M_{4}M_{8} +M_{4}M_{8}^{2},
\]
\begin{eqnarray*}
M_{15} 
& =  &  M_{8} +M_{4}^{2} -M_{4}M_{8} -M_{8}M_{4} +M_{8}^{2} -2M_{4}^{3} +3M_{4}^{2}M_{8}\\
& & \qquad {}-M_{4}M_{8}M_{4} +M_{8}^{2}M_{4} -M_{4}^{4} +M_{4}^{3}M_{8} +2M_{4}^{2}M_{8}M_{4} -2M_{4}^{2}M_{8}^{2},
\end{eqnarray*}
\[
M_{16} =  M_{4}M_{21}, \quad
M_{17} =  -I -M_{8} +M_{8}^{2}, \quad
M_{18} =  M_{8}M_{21},
\]
\[
M_{19} =  2M_{4}^{2} +M_{4}^{2}M_{8} +M_{4}M_{21}M_{8} -M_{4}^{2}M_{8}^{2},
\]
\[
M_{20} =  2M_{8}M_{4} +M_{4}^{3} -M_{4}^{2}M_{8} +M_{4}^{2}M_{21} -M_{4}M_{21}M_{8} -M_{4}^{3}M_{8} +M_{4}^{2}M_{8}^{2},
\]
\[
M_{22} =  M_{4} +M_{8} -M_{21} +M_{4}^{2} -M_{4}M_{21} -M_{4}M_{8}M_{4} +M_{4}M_{8}M_{21},
\]
\begin{eqnarray*}
M_{23} 
& = &
  I -M_{4} +M_{21} +M_{4}M_{8} +M_{8}M_{4} -M_{8}M_{21} -M_{21}M_{8} +M_{4}^{3} +M_{4}^{2}M_{8} \\
& & \qquad {}-M_{4}^{2}M_{21} +M_{4}M_{8}M_{4} -M_{4}M_{8}M_{21} -M_{4}M_{21}M_{8} -M_{4}^{4} +M_{4}^{3}M_{21},
\end{eqnarray*}
and
\[
M_{24} =  2M_{4} +M_{4}M_{8} +M_{21}M_{8} -M_{4}M_{8}^{2}.
\]
As before, a computer search was used to verify the cases
of Theorem~\ref{permisothm} corresponding to left cells with 
$\complexes W$-modules isomorphic to $M(\Gamma)=M(A_{19})$.
Moreover, the following holds.

\begin{theorem}
Suppose $\Gamma = A_{19}$  
and
$X$ is a nonempty subset of $\Gamma\cap \Gamma^{-1}$.
Then 
$\sum_{x \in X}\integers t_x$ is a subring (with 1) of
$J(\Gamma)$ if and only if
$X$ is one of sets
\[
\left\{x_2\right\}, \quad
\left\{x_1,x_{2}\right\}, \quad
\left\{x_2,x_{23}\right\}, \quad
\left\{x_1,x_2,x_{23},x_{24}\right\}, \quad
\left\{x_2,x_8,x_{17},x_{22},x_{23}\right\}, \quad
\Gamma\cap \Gamma^{-1}.
\]
\end{theorem}

Table~\ref{a19polys} contains the characteristic polynomials
of the matrices $M_j$.

\begin{table}
\caption{Characteristic polynomials for $\Gamma=A_{19}$.}
\[
\vbox{
\offinterlineskip
\halign{
\strut $#$ & \quad $#$ \hfil \cr
j & \det(u I - M_j) \cr
\noalign{\hrule}1 & (-1-u+u^{2})^{12}\cr
2 & (-1+u)^{24}\cr
3 & (1+u)^{4}(1-7u+u^{2})(-1-4u+u^{2})(1-3u+u^{2})^{4}(-1-u+u^{2})^{2}(-1+u+u^{2})(1+3u+u^{2})\cr
4 & (-1+u)^{2}(1+u)^{4}(1+u^{2})^{2}(-1-4u+u^{2})(1-3u+u^{2})(-1+u+u^{2})(1+3u^{2}+u^{4})^{2}\cr
5 & (-1-11u+u^{2})(-1-4u+u^{2})^{3}(-1-u+u^{2})^{3}(-1+u+u^{2})^{5}\cr
6 & (-1+u)^{2}(1+u)^{4}(1+u^{2})^{2}(-1-4u+u^{2})(1-3u+u^{2})(-1+u+u^{2})(1+3u^{2}+u^{4})^{2}\cr
7 & (-1+u)^{6}(1+u)^{4}(1-18u+u^{2})(1-7u+u^{2})(1-3u+u^{2})(1+3u+u^{2})^{4}\cr
8 & (1+u)^{6}(1-3u+u^{2})^{5}(-1-u+u^{2})^{4}\cr
9 & (-2+u)^{6}u^{8}(1+u)^{2}(2+u)^{2}(-4-8u+u^{2})(-1+4u+u^{2})^{2}\cr
10 & u^{8}(4-14u+u^{2})(-1+u+u^{2})(4+2u+u^{2})^{2}(-4+6u+7u^{2}-u^{3}+u^{4})^{2}\cr
11 & (-2+u)^{2}u^{8}(2+u)^{2}(-4-8u+u^{2})(4-6u+u^{2})^{2}(1-3u+u^{2})(-1+u+u^{2})^{2}\cr
12 & u^{8}(4-14u+u^{2})(-1+u+u^{2})(4+2u+u^{2})^{2}(-4+6u+7u^{2}-u^{3}+u^{4})^{2}\cr
13 & u^{8}(-4-22u+u^{2})(1-3u+u^{2})^{2}(-4+2u+u^{2})^{4}(1+3u+u^{2})\cr
14 & u^{8}(4-6u+u^{2})(-1-u+u^{2})(4+2u+u^{2})^{2}(-4-4u-3u^{2}-u^{3}+u^{4})^{2}\cr
15 & u^{8}(4-6u+u^{2})(-1-u+u^{2})(4+2u+u^{2})^{2}(-4-4u-3u^{2}-u^{3}+u^{4})^{2}\cr
16 & (-1+u)^{8}(1+u)^{2}(1-18u+u^{2})(1-3u+u^{2})(1+3u+u^{2})^{4}(1+7u+u^{2})\cr
17 & (-1+u)^{6}u^{8}(-4-2u+u^{2})^{5}\cr
18 & (-1-11u+u^{2})(-1-4u+u^{2})^{2}(-1-u+u^{2})^{2}(-1+u+u^{2})^{6}(-1+4u+u^{2})\cr
19 & (1+u)^{4}(1-7u+u^{2})(1-3u+u^{2})^{4}(-1-u+u^{2})(-1+u+u^{2})^{2}(1+3u+u^{2})(-1+4u+u^{2})\cr
20 & (-1+u)^{4}(1+u)^{2}(1+u^{2})^{2}(-1-4u+u^{2})(-1+u+u^{2})(1+3u+u^{2})(1+3u^{2}+u^{4})^{2}\cr
21 & (-1+u)^{4}(1+u)^{2}(1+u^{2})^{2}(-1-4u+u^{2})(-1+u+u^{2})(1+3u+u^{2})(1+3u^{2}+u^{4})^{2}\cr
22 & (1+u)^{6}(1-3u+u^{2})^{5}(-1+u+u^{2})^{4}\cr
23 & (-1+u)^{16}(1+u)^{8}\cr
24 & (-1-u+u^{2})^{8}(-1+u+u^{2})^{4}\cr
}}
\]
\label{a19polys}
\end{table}

%%%%%%%%%%%%%%%%%%%%%%%%%%%%%%%%%%%%%%%%%%%%%%%%%%%%%%
%%%%%%%%%%%%%%%%%%%%%%%%%%%%%%%%%%%%%%%%%%%%%%%%%%%%%%
%%   Section 7

\section{Concluding remarks}

We consider the effect of extending scalars to
 $K = \rationals[\sqrt{5}]$.
It is known that $K$ is a splitting
field for $W$ and $K(\sqrt{q})$ is a splitting
field for $\Hecke$ by \cite{AlvisLusztig}.  
For $\Gamma$ a left cell of $W$ and $F$ a field, put
\[
J(\Gamma)_{F} = F \otimes_\integers J(\Gamma).
\]
Since the coefficients of the structure constants
$h_{x,y,z}$ are nonnegative for $H_4$ 
by the calculation of du~Cloux \cite{Cloux}, a
result of Lusztig (\cite[21.9]{LusztigTwoParam})
shows that 
$J(\Gamma)_\complexes = \complexes \otimes_\integers J(\Gamma)$
is semisimple.
Thus $J(\Gamma)_{K}$ is semisimple.
A CAS program was used to compute the dimension of 
the derived algebra 
\[
\left[J(\Gamma)_\rationals,J(\Gamma)_\rationals\right]
=
\left< r s - s r \mid r, s \in J(\Gamma)_\rationals\right>.
\]
This dimension is 0 unless $\Gamma \subset A$, and is
$3$, $6$, and $12$ if
$\Gamma=A_1$, $\Gamma=A_9$, and $\Gamma=A_{19}$,
respectively.  Another CAS program has
verified that the number of central idempotents
in $J(\Gamma)_K$ is 11, 12, and 12 if
$\Gamma=A_1$, $\Gamma=A_9$, and $\Gamma=A_{19}$,
respectively.
From these observations and the structure
of the modules $K \otimes_\rationals M(\Gamma)$
given in \cite{AlvisHFour}, the following holds.

\begin{theorem}
Let $\Gamma$ be a left cell of $W$.
Then $J(\Gamma)_{K} = K \otimes_\integers J(\Gamma)$
is spit semisimple over $K$, and is isomorphic
to the endomorphism algebra of the $KW$-module
$K \otimes_\rationals M(\Gamma)$.  
\end{theorem}

Databases containing the
structure constants for $J(\Gamma)$
for all left cells $\Gamma$ and the
Kazhdan-Luztig polynomials are available 
from the author on request.

%%%%%%%%%%%%%%%%%%%%%%%%%%%%%%%%%%%%%%%%%%%%%%%%%%%
%%%%%   References

\end{document}